\documentclass[10pt,twoside,reqno]{amsart}
\usepackage{amssymb}
\usepackage{multirow}
\calclayout

\numberwithin{equation}{section}
\makeatletter

\renewcommand{\@secnumfont}{\bfseries}

\renewcommand{\section}{\@startsection{section}{1}%
  {0mm}{.7\linespacing\@plus\linespacing}{.5\linespacing}
  {\normalfont\bfseries\centering}}

\newcommand{\bibsection}{\@startsection{section}{1}%
  {0mm}{.7\linespacing\@plus\linespacing}{.5\linespacing}
  {\normalfont\scshape\centering}}

\renewcommand{\@biblabel}[1]{#1.}

\newtheorem{thm}{\bf Theorem}[section]

\begin{document}

\vspace{1.3cm}

\title {Degenerate Cauchy numbers and polynomials of the second kind}

\author{Taekyun Kim}
\address{Department of Mathematics, Kwangwoon University, Seoul 139-701, Republic
	of Korea}
\email{tkkim@kw.ac.kr}

\subjclass[2010]{11B68; 11S80}

\subjclass[2010]{11B68; 11S80}
\keywords{degenerate Cauchy numbers and polynomials}
\begin{abstract} 
Recently, degenerate Cauchy numbers and polynomials are introduced in [10]. In this paper, we study the degenerate Cauchy numbers and polynomials which are different from the previous degenerate Cauchy numbers and polynomials. In addition, we give some explicit identities for these numbers and polynomials which are derived from the generating function.
\end{abstract}
\maketitle
\bigskip
\medskip
\section{Introduction}
As is well known, the Cauchy polynomials are defined by the generating function to be
\begin{equation}\begin{split}\label{01}
\frac{t}{\log(1+t)}(1+t)^x = \sum_{n=0}^\infty C_n(x) \frac{t^n}{n!},\quad (\textnormal{see} \,\, [3,6,8,13]).
\end{split}\end{equation}

When $x=0$, $C_n=C_n(0)$ are called the Cauchy numbers. The higher-order Bernoulli polynomials are given by the generating function to be
\begin{equation}\begin{split}\label{02}
\left( \frac{t}{e^t-1} \right)^r e^{xt} = \sum_{n=0}^\infty B_n^{(r)} (x) \frac{t^n}{n!},\quad (\textnormal{see} \,\, [7,11,12]).
\end{split}\end{equation}
When $x=0$, $B_n^{(r)}=B_n^{(r)}(0)$ are called the higher-order Bernoulli numbers and $B_n^{(1)}(x) = B_n(x)$ are called the ordinary Bernoulli polynomials.

From \eqref{01} and \eqref{02}, we note that 
\begin{equation}\begin{split}\label{03}
C_n(x) = B_n^{(n)}(x+1),\,\,(n \geq 0),\quad (\textnormal{see} \,\, [10]).
\end{split}\end{equation}

In [1,2], L. Carlitz introduced the degenerate Bernoulli polynomials which are given by the generating function to be

\begin{equation}\begin{split}\label{04}
\frac{t}{(1+\lambda t)^{\frac{1}{\lambda }}-1}(1+\lambda t)^{\frac{x}{\lambda }} = \sum_{n=0}^\infty \beta_{n,\lambda }(x) \frac{t^n}{n!}.
\end{split}\end{equation}

When $x=0$, $\beta_{n,\lambda }=\beta_{n,\lambda }(0)$ are called the degenerate Bernoulli numbers.

Note that $\lim_{\lambda  \rightarrow 0} \beta_{n,\lambda }(x) = B_n(x)$, $(n \geq 0)$. In the viewpoint of \eqref{04}, the degenerate Cauchy polynomials are defined by the generating function to be

\begin{equation}\begin{split}\label{05}
\int_0^1 (1+ \log(1+\lambda t)^{\frac{1}{\lambda }})^{x+y} dy&= \frac{\log(1+\lambda t)^{\frac{1}{\lambda }}}{\log(1+\log(1+\lambda t)^{\frac{1}{\lambda }})}(1+\log(1+\lambda t)^{\frac{1}{\lambda }})^x\\
&=\sum_{n=0}^\infty C_{n,\lambda }^* (x) \frac{t^n}{n!},\quad (\textnormal{see} \,\, [10]).
\end{split}\end{equation}

The falling factorial sequences are given by 

\begin{equation}\begin{split}\label{06}
(x)_0=1, (x)_n = x(x-1)\cdots(x-n+1),\,\,(n \geq 1),\quad (\textnormal{see} \,\, [4,5,6]).
\end{split}\end{equation}
The Stirling number of the first kind is defined as
\begin{equation}\begin{split}\label{07}
(x)_n = \sum_{l=0}^n S_1(n,l) x^l,\,\,(n \geq 0),\quad (\textnormal{see} \,\, [8,9,10]),
\end{split}\end{equation}
and the Stirling number of the second kind is given by 
\begin{equation}\begin{split}\label{08}
x^n = \sum_{l=0}^n S_2(n,l)(x)_l, \,\,(n \geq 0).
\end{split}\end{equation}

In [10], it were known that
\begin{equation}\begin{split}\label{09}
C_m(x) = \sum_{n=0}^m C_{m,\lambda }^* \lambda ^{m-n} S_2(m,n),
\end{split}\end{equation}
and

\begin{equation}\begin{split}\label{10}
C_{m,\lambda }^* (x) = \sum_{n=0}^m B_n^{(n)} (x+1) \lambda ^{m-n} S_1(m,n).
\end{split}\end{equation}

In this paper, we study the degenerate Cauchy numbers and polynomials of the second kind which are different from previous degenerate Cauchy numbers and polynomials. In addition, we give some explicit identities for the degenerate Cauchy numbers and polynomials of the second kind which are derived from the generating function.

\section{Degenerate Cauchy numbers and polynomials of the second kind}

Now, we define the degenerate Cauchy polynomials of the second kind which are given by the generating function to be

\begin{equation}\begin{split}\label{11}
\frac{t}{\log \big(1+ \frac{1}{\lambda } \log(1+\lambda t)\big)} \big(1+\tfrac{1}{\lambda }\log(1+\lambda  t)\big)^x = \sum_{n=0}^\infty C_{n,\lambda }(x) \frac{t^n}{n!}.
\end{split}\end{equation}

When $x=0$, $C_{n,\lambda }= C_{n,\lambda }(0)$ are called the degenerate Cauchy numbers of the second kind.

By replacing $t$ by $\frac{1}{\lambda }(e^{\lambda t}-1)$ in \eqref{11}, we get

\begin{equation}\begin{split}\label{12}
&\sum_{m=0}^\infty C_{m,\lambda }(x) \lambda ^{-m} \frac{1}{m!} (e^{\lambda t}-1)^m = \left( \frac{1}{\lambda } \frac{t}{\log(1+t)} (1+t)^x \right) \left(\frac{e^{\lambda t }-1}{t} \right)\\
&= \left( \frac{1}{\lambda } \sum_{m=0}^\infty C_m(x) \frac{t^m}{m!} \right) \left( \frac{1}{t} \sum_{l=1}^\infty \frac{\lambda ^l t^l}{l!} \right) =\left( \frac{1}{\lambda } \sum_{m=0}^\infty C_m(x) \frac{t^m}{m!} \right) \left( \frac{1}{t} \sum_{l=0}^\infty \frac{\lambda^{l+1}}{l+1} \frac{t^l}{l!} \right)\\
&=\sum_{n=0}^\infty \left( \sum_{m=0}^n {n \choose m} C_{n-m}(x) \frac{\lambda ^m}{m+1} \right) \frac{t^n}{n!}.
\end{split}\end{equation}

On the other hand,
\begin{equation}\begin{split}\label{13}
\sum_{m=0}^\infty C_{m,\lambda }(x) \lambda ^{-m} \sum_{n=m}^\infty S_2(n,m) \lambda ^n \frac{t^n}{n!} = \sum_{n=0}^\infty \left( \sum_{m=0}^n \lambda ^{n-m} C_{m,\lambda }(x) S_2(n,m) \right) \frac{t^n}{n!}.
\end{split}\end{equation}

Therefore, by \eqref{12} and \eqref{13}, we obtain the following theorem.

\begin{thm}
For $n \geq 0$, we have
\begin{equation*}\begin{split}
 \sum_{m=0}^n \lambda ^{n-m} C_{m,\lambda }(x) S_2(n,m) = \sum_{m=0}^n {n \choose m} C_{n-m}(x) \frac{\lambda ^m}{m+1}.
\end{split}\end{equation*}
\end{thm}

From Theorem 1, we note that $\lim_{\lambda  \rightarrow 0} C_{m,\lambda }(x) = C_{n-m}(x)$. Let us take $t= \frac{1}{\lambda } \log(1+\lambda t)$ in \eqref{01}. Then we have

\begin{equation}\begin{split}\label{14}
&\sum_{m=0}^\infty C_m(x) \lambda ^{-m} \frac{1}{m!} \big(\log(1+\lambda t)\big)^m = \frac{\frac{1}{\lambda }\log(1+\lambda t) }{\log\big(1+\frac{1}{\lambda }\log(1+\lambda t)\big)} \big(1+\tfrac{1}{\lambda } \log(1+\lambda t)\big)^x \\
&=\left( \frac{\log(1+\lambda t)}{\lambda t}\right)
 \left( \frac{t}{\log\big( 1+\frac{1}{\lambda }\log(1+\lambda t)\big) } \log\big( 1+\tfrac{1}{\lambda }\log(1+\lambda t)\big)^x	\right)\\
&=\left( \sum_{l=1}^\infty \frac{(-1)^{l-1}}{l} (\lambda t)^{l-1} \right) \left( \sum_{m=0}^\infty C_{m,\lambda }(x) \frac{t^m}{m!} \right) \\
&= \left( \sum_{l=0}^\infty \frac{\lambda ^l(-1)^l}{l+1}t^l \right) \left( \sum_{m=0}^\infty C_{m,\lambda }(x)\frac{t^m}{m!} \right)\\
&= \sum_{n=0}^\infty \left( \sum_{m=0}^n \frac{n! \lambda ^{n-m} (-1)^{n-m}}{(n-m+1) m!} C_{m,\lambda }(x) \right) \frac{t^n}{n!} \\
&= \sum_{n=0}^\infty \left( \sum_{m=0}^n {n \choose m} \frac{(n-m)!}{n-m+1} \lambda ^{n-m} (-1)^{n-m} C_{m,\lambda }(x) \right) \frac{t^n}{n!}.
\end{split}\end{equation}

On the other hand

\begin{equation}\begin{split}\label{15}
\sum_{m=0}^\infty C_m(x) \lambda ^{-m} \big( \log(1+\lambda t) \big)^m &=
\sum_{m=0}^\infty C_m(x) \lambda ^{-m} \sum_{n=m}^\infty S_1(n,m) \frac{\lambda ^n}{n!} t^n\\
&= \sum_{n=0}^\infty \left( \sum_{m=0}^n C_m(x) \lambda ^{n-m} S_1(n,m) \right) \frac{t^n}{n!}.
\end{split}\end{equation}
Therefore, by \eqref{14} and \eqref{15}, we obtain the following theorem.

\begin{thm}
For $n \geq 0$, we have

\begin{equation*}\begin{split}
 \sum_{m=0}^n C_m(x) \lambda ^{n-m} S_1(n,m) = \sum_{m=0}^n  \frac{(n-m)!}{n-m+1}{n \choose m} (-1)^{n-m} \lambda ^{n-m}  C_{m,\lambda }(x).
\end{split}\end{equation*}
\end{thm}

When $x=0$, $C_m=C_m(0)$ and $C_{m,\lambda} = C_{m,\lambda }(0)$. So, by Theorem 2, we get
\begin{equation*}\begin{split}
 \sum_{m=0}^n C_m \lambda ^{n-m} S_1(n,m) = \sum_{m=0}^n  \frac{(n-m)!}{n-m+1}{n \choose m} (-1)^{n-m} \lambda ^{n-m}  C_{m,\lambda }.
\end{split}\end{equation*}

From \eqref{11}, we note that

\begin{equation}\begin{split}\label{16}
\sum_{n=0}^\infty C_{n,\lambda }(x) \frac{t^n}{n!} &= \left( \sum_{l=0}^\infty C_{l,\lambda } \frac{t^l}{l!} \right) \left( \sum_{m=0}^\infty (x)_m \lambda ^{-m} \frac{1}{m!} \big( \log(1+\lambda t) \big)^m \right)\\
&= \left( \sum_{l=0}^\infty C_{l,\lambda } \frac{t^l}{l!} \right)\left( \sum_{m=0}^\infty (x)_m \lambda ^{-m} \sum_{k=m}^\infty S_1(k,m) \lambda ^k \frac{t^k}{k!} \right)
\\&= \left( \sum_{l=0}^\infty C_{l,\lambda } \frac{t^l}{l!} \right)\left( \sum_{k=0}^\infty \left( \sum_{m=0}^k (x)_m \lambda ^{k-m} S_1(k,m) \right) \frac{t^k}{k!} \right)\\
&= \sum_{n=0}^\infty \left( \sum_{k=0}^n \sum_{m=0}^k {n \choose k} C_{n-k} (x)_m \lambda ^{k-m} S_1(k,m) \right) \frac{t^n}{n!}.
\end{split}\end{equation}

Therefore, by comparing the coefficients on the both sides of \eqref{16}, we obtain the following theorem.

\begin{thm}
For $n \geq 0$, we have
\begin{equation*}\begin{split}
C_{n,\lambda }(x) =  \sum_{k=0}^n \sum_{m=0}^k {n \choose k}(x)_m C_{n-k}  \lambda ^{k-m} S_1(k,m).
\end{split}\end{equation*}
\end{thm}

From \eqref{11}, we note that
\begin{equation}\begin{split}\label{17}
t &= \left( \sum_{l=0}^\infty C_{l,\lambda } \frac{t^l}{l!} \right) \left( \log \big(1+ \tfrac{1}{\lambda } \log(1+\lambda t)\big) \right)\\
 &= \left( \sum_{l=0}^\infty C_{l,\lambda } \frac{t^l}{l!} \right) \left( \sum_{k=1}^\infty \frac{(-1)^{k-1}}{k} \lambda ^{-k} \big( \log(1+\lambda t) \big)^k \right)\\
 &= \left( \sum_{l=0}^\infty C_{l,\lambda } \frac{t^l}{l!} \right) \left( \sum_{k=1}^\infty (k-1)! (-1)^{k-1} \lambda ^{-k} \sum_{m=k}^\infty S_1(m,k) \lambda ^m \frac{t^m}{m!} \right)
 \\ &= \left( \sum_{l=0}^\infty C_{l,\lambda } \frac{t^l}{l!} \right) \left( \sum_{m=1}^\infty \left( \sum_{k=1}^m (k-1)! (-1)^{k-1} \lambda ^{m-k} S_1(m,k) \right) \frac{t^m}{m!} \right)\\
 &= \sum_{n=1}^\infty \left( \sum_{m=1}^n \sum_{k=1}^m {n \choose m} C_{n-m,\lambda } (k-1)! (-1)^{k-1} \lambda ^{m-k} S_1(m,k) \right) \frac{t^n}{n!}.
\end{split}\end{equation}
Comparing the coefficients on the both sides of \eqref{17}, we obtain the following theorem.

\begin{thm}
For $n \in \mathbb{N}$, we have
\begin{equation*}\begin{split}
C_{0,\lambda }=1, \,\,\sum_{m=1}^n \sum_{k=1}^m {n \choose m} C_{n-m,\lambda } (k-1)! (-1)^{k-1} \lambda ^{m-k} S_1(m,k)= \begin{cases}
1,&\text{if}\,\,n=1,\\
0,&\text{if}\,\,n>1.
\end{cases}
\end{split}\end{equation*}
\end{thm}

Now, we observe that

\begin{equation}\begin{split}\label{18}
\frac{t}{\log \big( 1+ \tfrac{1}{\lambda }\log(1+\lambda t) \big)} &= \left( \frac{\lambda t}{\log(1+\lambda t)}\right) \left( \frac{\frac{1}{\lambda }\log(1+\lambda t)}{\log \big( 1+ \tfrac{1}{\lambda }\log(1+\lambda t) \big)} \right) \\
&=\left( \sum_{l=0}^\infty \lambda ^l B_l^{(l)} \frac{t^l}{l!} \right) \left( \sum_{m=0}^\infty C_{m,\lambda }^* \frac{t^m}{m!} \right)\\
&= \sum_{n=0}^\infty \left( \sum_{m=0}^n {n \choose m} \lambda ^{n-m} B_{n-m}^{(n-m)} C_{m,\lambda }^* \right) \frac{t^n}{n!}.
\end{split}\end{equation}
From \eqref{11}, we note that
\begin{equation}\begin{split}\label{19}
\frac{t}{\log \big( 1+ \tfrac{1}{\lambda }\log(1+\lambda t) \big)} = \sum_{n=0}^\infty C_{n,\lambda } \frac{t^n}{n!}.
\end{split}\end{equation}

Therefore, by \eqref{18} and \eqref{19}, we obtain the following theorem.

\begin{thm}
For $n \geq 0$, we have
\begin{equation*}\begin{split}
C_{n,\lambda }= \sum_{m=0}^n {n \choose m} \lambda ^{n-m} B_{n-m}^{(n-m)} C_{m,\lambda }^*.
\end{split}\end{equation*}
\end{thm}
By \eqref{05}, we get

\begin{equation}\begin{split}\label{20}
\sum_{n=0}^\infty C_{n,\lambda }^* \frac{t^n}{n!} &= \frac{\frac{1}{\lambda }\log(1+\lambda t) }{\log \big( 1+ \tfrac{1}{\lambda }\log(1+\lambda t) \big)} = \left( \frac{\log(1+\lambda t)}{\lambda t} \right) \left( \frac{t}{\log \big( 1+ \tfrac{1}{\lambda }\log(1+\lambda t) \big)} \right)\\
&= \left( \sum_{l=0}^\infty \lambda ^l D_l \frac{t^l}{l!} \right) \left( \sum_{m=0}^\infty C_{m,\lambda } \frac{t^m}{m!} \right) \\
&= \sum_{n=0}^\infty \left( \sum_{m=0}^n {n \choose m} \lambda ^{n-m} D_{n-m} C_{m,\lambda } \right) \frac{t^n}{n!},
\end{split}\end{equation}
Where $D_n$, $(n \geq 0)$, are the Daehee numbers which are defined by the generating function to be 

\begin{equation*}\begin{split}
\frac{\log(1+t)}{t} = \sum_{n=0}^\infty D_n \frac{t^n}{n!},\quad (\textnormal{see} \,\, [4,5,9,11,12]).
\end{split}\end{equation*}

Comparing the coefficients on the both sides of \eqref{20}, we obtain the following theorem.

\begin{thm}
For $n \geq 0$, we have
\begin{equation*}\begin{split}
C_{n,\lambda }^* = \sum_{m=0}^n {n \choose m} \lambda ^{n-m} D_{n-m} C_{m,\lambda }.
\end{split}\end{equation*}
\end{thm}

From \eqref{11}, we note that

\begin{equation}\begin{split}\label{21}
\sum_{n=0}^\infty C_{n,\lambda }(1) \frac{t^n}{n!} &= \frac{t}{\log \big( 1+ \tfrac{1}{\lambda }\log(1+\lambda t) \big)} \big( 1+ \tfrac{1}{\lambda }\log(1+\lambda t) \big)\\
&= \sum_{n=0}^\infty C_{n,\lambda } \frac{t^n}{n!} + t \cdot \frac{\frac{1}{\lambda }\log(1+\lambda t)}{\log \big( 1+ \tfrac{1}{\lambda }\log(1+\lambda t) \big)}.
\end{split}\end{equation}
By \eqref{05} and \eqref{21}, we get

\begin{equation}\begin{split}\label{22}
t \sum_{n=0}^\infty C_{n,\lambda }^* \frac{t^n}{n!} &= t \frac{\frac{1}{\lambda }\log(1+\lambda t)}{\log \big( 1+ \tfrac{1}{\lambda }\log(1+\lambda t) \big)} = \sum_{n=1}^\infty \Big( C_{n,\lambda }(1) - C_{n,\lambda } \Big) \frac{t^n}{n!}\\
&=\sum_{n=0}^\infty \left( \frac{C_{n+1,\lambda }(1) - C_{n+1,\lambda }}{n+1} \right) \frac{t^{n+1}}{n!}.
\end{split}\end{equation}

From \eqref{22}, we get
\begin{equation}\begin{split}\label{23}
\sum_{n=0}^\infty C_{n,\lambda }^* \frac{t^n}{n!} = \sum_{n=0}^\infty \left( \frac{C_{n+1,\lambda }(1) - C_{n+1,\lambda }}{n+1} \right) \frac{t^n}{n!}.
\end{split}\end{equation}

Comparing the coefficients on the both sides of \eqref{23}, we obtain the following theorem.

\begin{thm}
For $n \geq 0$, we have
\begin{equation*}\begin{split}
C_{n,\lambda }^* =  \frac{C_{n+1,\lambda }(1) - C_{n+1,\lambda }}{n+1}.
\end{split}\end{equation*}

\end{thm}

By \eqref{05}, we get

\begin{equation}\begin{split}\label{24}
\sum_{n=0}^\infty C_{n,\lambda }^* (1) \frac{t^n}{n!} &= \frac{\frac{1}{\lambda}\log(1+\lambda t)}{\log \big( 1+ \tfrac{1}{\lambda }\log(1+\lambda t) \big)} \big( 1+ \tfrac{1}{\lambda }\log(1+\lambda t) \big)\\
&= \sum_{n=0}^\infty C_{n,\lambda }^* \frac{t^n}{n!} + \left( \frac{t^2}{\log \big( 1+ \tfrac{1}{\lambda }\log(1+\lambda t) \big)} \right) \left( \frac{\log(1+\lambda t)}{\lambda t} \right)^2.
\end{split}\end{equation}

For $r \in \mathbb{N}$, the higher-order Daehee numbers are defined by the generating function to be
\begin{equation}\begin{split}\label{25}
\left( \frac{\log(1+t)}{t} \right)^r  = \sum_{n=0}^\infty D_n^{(r)} \frac{t^n}{n!},\quad (\textnormal{see} \,\, [4,7]).
\end{split}\end{equation}

From \eqref{24} and \eqref{25}, we have

\begin{equation}\begin{split}\label{26}
\sum_{n=1}^\infty \Big( C_{n,\lambda }^* (1) - C_{n,\lambda }^* \Big) \frac{t^n}{n!} &=
t \left( \frac{t}{\log \big( 1+ \tfrac{1}{\lambda }\log(1+\lambda t) \big)} \right) \left( \frac{\log(1+\lambda t)}{\lambda t} \right)^2\\
&=t \left( \sum_{l=0}^\infty C_{l,\lambda } \frac{t^l}{l!} \right) \left( \sum_{m=0}^\infty D_m^{(2)} \frac{t^m}{m!} \right)\\
&= t \left( \sum_{n=0}^\infty \left( \sum_{l=0}^n {n \choose l} C_{l,\lambda } D_{n-l,\lambda }^{(2)} \right) \frac{t^n}{n!} \right).
\end{split}\end{equation}

Thus, by \eqref{26}, we get

\begin{equation}\begin{split}\label{27}
\sum_{n=0}^\infty \left( \frac{C_{n+1,\lambda }^*(1)-C_{n+1,\lambda }^*}{n+1} \right) \frac{t^n}{n!} =  \sum_{n=0}^\infty \left( \sum_{l=0}^n {n \choose l} C_{l,\lambda } D_{n-l,\lambda }^{(2)} \right) \frac{t^n}{n!}.
\end{split}\end{equation} 
Comparing the coefficients on the both sides of \eqref{27}, we obtain the following theorem.

\begin{thm}
For $n \geq 0$, we have
\begin{equation*}\begin{split}
 \frac{C_{n+1,\lambda }^*(1)-C_{n+1,\lambda }^*}{n+1} =  \sum_{l=0}^n {n \choose l} C_{l,\lambda } D_{n-l,\lambda }^{(2)}.
\end{split}\end{equation*}

\end{thm}

\end{document}